\documentclass[reqno,12pt]{amsart}

\usepackage{amsmath}
\usepackage{amssymb}
\usepackage{amsthm}
\usepackage[english]{babel}
\newtheorem{theorem}{Theorem}

\usepackage{fancyhdr}  
\newtheorem{lemma}{Lemma}

\newtheorem{corollary}{Corollary}

\begin{document}

\title {On the Littlewood--Offord
problem}

\author{Yulia S. Eliseeva$^{1,2}$, Andrei Yu. Zaitsev$^{1,3}$}


\email{pochta106@yandex.ru}
\address{St.~Petersburg State University\newline\indent
and Laboratory of
 Chebyshev in St. Petersburg State University \bigskip}

\email{zaitsev@pdmi.ras.ru}
\address{St.~Petersburg Department of Steklov Mathematical Institute
\newline\indent
Fontanka 27, St.~Petersburg 191023, Russia\newline\indent and
St.~Petersburg State University}

\begin{abstract}{The paper deals with studying a connection of the
Littlewood--Offord problem with estimating the concentration
functions of some symmetric infinitely divisible distributions.
Some multivariate generalizations of results of Arak (1980) are
given. They show a connection of the concentration function of the
sum with the arithmetic structure of supports of distributions of
independent random vectors for arbitrary distributions of
summands.}\end{abstract}


\maketitle

\footnotetext[1]{The first and the second authors are supported
by grants RFBR 13-01-00256 and NSh-2504.2014.1.}
\footnotetext[2]{The first author was supported by Laboratory of
 Chebyshev in St. Petersburg State University (grant of the Government of Russian Federation 11.G34.31.0026)
 and
by grant of St. Petersburg State University
6.38.672.2013.}\footnotetext[3]{The second author was supported by
the Program of Fundamental Researches of Russian Academy of
Sciences "Modern Problems of Fundamental Mathematics".}

Let $X,X_1,\ldots,X_n$ be independent identically distributed
(i.i.d.) random variables. Let $a=(a_1,\ldots,a_n)$, where
$a_k=(a_{k1},\ldots,a_{kd})\in \mathbf{R}^d$, $k=1,\ldots, n$. The
 concentration function of a $\mathbf{R}^d$-dimensional random
vector $Y$ with distribution $F=\mathcal L(Y)$ is defined by the
equality
\begin{equation}
Q(F,\lambda)=\sup_{x\in\mathbf{R}^d}\mathbf{P}(Y\in x+ \lambda B),
\quad \lambda\geq0, \nonumber
\end{equation}
where $B=\{x\in\mathbf{R}^n:\|x\|\leq 1/2\}$. In this paper we
study the behavior of the concentration functions of the weighted
sums $S_a=\sum\limits_{k=1}^{n} X_k a_k$ with respect to the
properties of vectors~$a_k$. Recently, interest in this subject
has increased considerably in connection with the study of
eigenvalues of random matrices (see, for instance, Friedland and
Sodin \cite{Fried:Sod:2007}, Nguyen and Vu~\cite{Nguyen:Vu:2011},
Rudelson and Vershynin \cite{Rud:Ver:2008}, \cite{Rud:Ver:2009},
Tao and Vu \cite{Tao:Vu:2009:Ann}, \cite{Tao:Vu:2009:Bull},
 Vershynin \cite{Ver:2011}). For a detailed history of
the problem we refer to a recent review of Nguen and Vu
\cite{Nguyen and Vu13}. The authors of the above articles (see
also Hal\'asz \cite{Hal:1977}) called this question the
Littlewood--Offord problem, since, for the first time, this
problem was considered in 1943 by Littlewood and Offord
\cite{Lit:Off:1943}  in connection with the study of random
polynomials. They considered a special case, where the
coefficients $a_k \in \mathbf{R}$ are one-dimensional, and $X$
takes values $\pm1$ with probabilities~$1/2$.

Let us introduce some notation. In the sequel, let $F_a$ denote
the distribution of the sum $S_a$, $E_y$ is the probability
measure concentrated at a point $y$, and let $G$ be the
distribution of the random variable $\widetilde{X}$, where
$\widetilde{X}=X_1-X_2$ is the symmetrized random variable.

The symbol $c$ will be used for absolute positive constants which
may be different even in the same formulas.

Writing $A\ll B$ means that $|A|\leq c B$. Also we will write
$A\asymp B$, if $A\ll B$ and $B\ll A$. We will write $A\ll _{d}
B$, if $|A|\leq c(d) B$, where $c(d)>0$ depends on $d$ only.
Similarly, $A\asymp_{d} B$, if $A\ll_{d} B$ and $B\ll_{d} A$. The
scalar product in $\mathbf{R}^d$ will be denoted $\left\langle \,
\cdot\,,\,\cdot\,\right\rangle$. Later $\lfloor x\rfloor$ is the
largest integer~$k$ such that~$k< x$.
 For~${x=(x_1,\dots,x_n )\in\mathbf R^n}$ we will use the norms
 $\|x\|^2= x_1^2+\dots +x_n^2$ and $|x|=
\max_j|x_j|$. We denote by $\widehat F(t)$, $t\in\mathbf R^d$, the
characteristic function of $d$-dimensional distributions~$F$.

Products and powers of measures will be understood in the
convolution sense. While a distribution~$F$ is infinitely
divisible, $F^\lambda$, $\lambda\ge0$, is the infinitely divisible
distribution with characteristic function $\widehat F^\lambda(t)$.

The elementary properties of concentration functions are well
studied (see, for instance, \cite{Arak:Zaitsev:1986},
\cite{Hen:Teod:1980}, \cite{Petrov:1972}). It is known that
\begin{equation}\label{8jrt}
Q(F,\mu)\ll_d (1+ \lfloor \mu/\lambda \rfloor)^d\,Q(F,\lambda)
\end{equation}
 for any $\mu,\lambda>0$.
Hence,
\begin{equation}\label{8art}
Q(F,c\lambda)\asymp_{d}\,Q(F,\lambda).
\end{equation}

Let us formulate a generalization of the classical Ess\'een
inequality \cite{Ess:1966} to the multivariate case
(\cite{Ess:1968}, see also \cite{Hen:Teod:1980}):

\begin{lemma}\label{lm3} Let $\tau>0$ and let
 $F$ be a $d$-dimensional probability distribution. Then
\begin{equation}
Q(F, \tau)\ll_d \tau^d\int_{|t|\le1/\tau}|\widehat{F}(t)| \,dt.
\label{4s4d}
\end{equation}
\end{lemma}

In the general case $Q(F,\lambda)$ cannot be estimated from below
by the right hand side of
 inequality~\eqref{4s4d}. However, if we assume additionally that the distribution $F$ is symmetric and
 its characterictic function is non-negative for all~$t\in\mathbf R$, then we have the lower bound:
\begin{equation} \label{1a}Q(F, \tau)\gg_d \tau^d\int\limits_{|t|\le1/\tau}{|\widehat{F}(t)|
\,dt},
\end{equation}
and, therefore,
\begin{equation} \label{1b}
Q(F, \tau)\asymp_d
\tau^d\int\limits_{|t|\le1/\tau}{|\widehat{F}(t)| \,dt},
\end{equation} (see \cite{Arak:1980} or \cite{Arak:Zaitsev:1986}, Lemma~1.5 of Chapter II for $d=1$).
In the multivariate case relations \eqref{1a} and~\eqref{1b} were
obtained by Zaitsev \cite{Zaitsev:1987}, see also Eliseeva
\cite{Eliseeva}. Just the use of relation \eqref{1b} allows us to
simplify the arguments of Friedland and Sodin
 \cite{Fried:Sod:2007},  Rudelson and Vershynin~\cite{Rud:Ver:2009} and Vershynin \cite{Ver:2011}
 which were applied to  Littlewood--Offord problem
 (see \cite{Eliseeva}, \cite{EGZ} and \cite{Eliseeva and Zaitsev}).
 \medskip

The main result of this paper is a general inequality which
reduces the estimation of concentration functions in the
Littlewood--Offord problem to the estimation of concentration
functions of some infinitely divisible distributions. This result
is formulated in Theorem~\ref{lm43}.

For $z\in \mathbf{R}$, introduce the distribution $H_z$ with the
characteristic function
\begin{equation} \label{11b}\widehat{H}_z(t)
=\exp\Big(-\frac{\,1\,}2\;\sum_{k=1}^{n}\big(1-\cos(\left\langle
\, t,a_k\right\rangle z)\big)\Big).
\end{equation}It depends on the vector~$a$.
It is clear that $H_z$ is a symmetric infinitely divisible
distribution.
  Therefore, its characteristic function is positive for all $t\in \mathbf{R}^d$.

\begin{theorem}\label{lm43}Let\/ $V$ be an arbitrary\/
$d$-dimensional Borel measure such that\/ $\lambda=V\{\mathbf
R\}>0$,\/ and $V\le G$, that is,\/ $V\{B\}\le G\{B\}$, for any
Borel set~$B$. Then, for any\/ $\varepsilon>0$ and\/ $\tau>0$, we
have
\begin{equation}\label{cc23}
Q(F_a,\tau) \ll_d Q(H_{1}^{\lambda},\varepsilon)\,\exp\bigg(d
\int_{z\in\mathbf{R}}\log\big(1+\lfloor\tau(\varepsilon
|z|)^{-1}\rfloor \big)\,F\{dz\}\bigg),
\end{equation}where $F=\lambda^{-1}V$.\end{theorem}

Note that $\log\big(1+\lfloor\tau(\varepsilon |z|)^{-1}\rfloor
\big)=0$, for $|z|\ge\tau/\varepsilon$. Therefore, the integration
in~\eqref{cc23} is taken, in fact, over the set
$\big\{z:|z|<\tau/\varepsilon\big\}$ only.

\begin{corollary}\label{lm42}Let\/ $\delta>0$ and
\begin{equation}p(\delta)=
G\big\{\{z:|z| \ge \delta\}\big\}>0.\end{equation}Then, for any\/
$\varepsilon,\tau>0$, we have
\begin{equation}\label{10abc}
Q(F_a, \tau) \ll_d e^{\Delta}\, Q(H_1^{p(\delta)}, \varepsilon),
\end{equation}
where
\begin{equation}\Delta=\Delta(\tau,\varepsilon,\delta)=\frac{d}{p(\delta)}\int\limits_{|z|
\ge \delta} \log\big(1+\lfloor\tau(\varepsilon |z|)^{-1}\rfloor
\big) \, G\{dz\}.\end{equation}
\end{corollary}

In particular, choosing $\delta=\tau/\varepsilon$, we get
\begin{corollary}\label{lm452}For any\/ $\varepsilon,\tau>0$,
we have
\begin{equation}\label{11abc}
Q(F_a, \tau) \ll_d Q(H_1^{p(\tau/\varepsilon)}, \varepsilon).
\end{equation}
\end{corollary}Just the statement of Corollary \ref{lm452} (usually for $\tau=\varepsilon$)
 is actually the starting point of almost
all recent studies on the Littlewood--Offord problem (see, for
instance, \cite{Fried:Sod:2007},~\cite{Hal:1977},
\cite{Nguyen:Vu:2011}, \cite{Rud:Ver:2008}, \cite{Rud:Ver:2009}
and \cite{Ver:2011}). More precisely, with the help of Lemma
\ref{lm3} or its analogs, the authors of the above-mentioned
papers have obtained estimates of the type
\begin{equation}\label{11abc2}
Q(F_a, \tau) \ll_d \sup_{z\ge\tau/\varepsilon
}\tau^d\int_{|t|\le1/\tau}\widehat{H}_z^{p(\tau/\varepsilon)}(t)
\,dt.
\end{equation}
The fact that  \eqref{8jrt} and \eqref{1b} imply that
\begin{eqnarray}
\sup_{z\ge\tau/\varepsilon
}\tau^d\int_{|t|\le1/\tau}\widehat{H}_z^{p(\tau/\varepsilon)}(t)
\,dt  &\asymp_d&\sup_{z\ge\tau/\varepsilon}\; Q\big(H_{z}^{p(\tau/\varepsilon)},\tau\big)\nonumber \\
= \sup_{z\ge\tau/\varepsilon}\;
Q\big(H_{1}^{p(\tau/\varepsilon)},\tau/z\big)&=&Q\big(H_{1}^{p(\tau/\varepsilon)},\varepsilon\big),
\end{eqnarray}
 remained apparently
unnoticed by the authors of these papers that significantly
hampered further evaluation of the right-hand side of
inequality~\eqref{11abc2}.

Choosing $V$ so that
\begin{equation}\label{cc243}
V\{dz\}=\big(\max\big\{1,\,\log\big(1+\lfloor\tau(\varepsilon
|z|)^{-1}\rfloor \big)\big\}\big)^{-1}\,G\{dz\},
\end{equation}
we obtain
\begin{corollary}\label{lm429}
For any\/ $\varepsilon,\tau>0$, we have
\begin{equation}\label{cc234}
Q(F_a,\tau) \ll_d
Q(H_{1}^{\lambda},\varepsilon)\,\exp\big(d\lambda^{-1}G\big\{\{z:|z|<\tau/\varepsilon\}\big\}
\big),
\end{equation}where
\begin{equation}\label{cc239}
\lambda=\lambda(G,\tau/\varepsilon)=V\{\mathbf
R\}=\int\limits_{z\in\mathbf{R}}\big(\max\big\{1,\,\log\big(1+\lfloor\tau(\varepsilon
|z|)^{-1}\rfloor \big)\big\}\big)^{-1}\,G\{dz\}.
\end{equation}
\end{corollary}

In Corollaries~\ref{lm42}--\ref{lm429} we choose the measure $V$
in the form $V\{dz\}=f(z)\,G\{dz\}$ with $0\le f(z)\le1$. It is
not clear what choice of $f$ is optimal. This depends on $a$ and
$G$.

Choosing the optimal function $f$, minimizing the right-hand sides
of inequalities \eqref{10abc}, \eqref{11abc} and \eqref{cc234}, is
a difficult problem. It is clear that its solution depends on $a$
and $G$. Certainly, it is sufficient to consider non-decreasing
functions~$f$ only.

For a fixed $\varepsilon$, an increase of~$\lambda$ implies a
decrease of $Q(H_{1}^{\lambda},\varepsilon)$. Theorem \ref{lm43}
may be applied for $V=G$. Then $\lambda=1$. This is the maximal
possible value of~$\lambda$. However, the integral in the
right-hand side of~\eqref{cc23} may be in this case infinite. In
particular, it diverges if the distribution~$G$ has a nonzero atom
at zero. This atom in any case should be excluded in constructing
the measure $V$, if we expect to get a meaningful bound for
 $Q(F_a,\tau)$. For a fixed measure~$V$, decreasing $\varepsilon$ implies
a decrease of $Q(H_{1}^{\lambda},\varepsilon)$, but an increase of
the integral in the right-hand side of inequality~\eqref{cc23}.

In Corollary \ref{lm429} we used the measure $V$, defined in
\eqref{cc243} so that the integral in the right-hand side of
inequality~\eqref{cc23} would converge always, no matter what is
the measure~$G$.

The proof of Theorem~$\ref{lm43}$ is based on elementary
properties of concentration functions, it will be given below.
Note that $H_1^{\lambda}$ is an infinitely divisible distribution
with the L\'evy
 spectral measure $M_\lambda=\frac{\,\lambda\,}4\;M^*$, where
 $M^*=\sum\limits_{k=1}^{n}\big(E_{a_k}+E_{-a_k}\big)$.
It is clear that the assertions of Theorem~$\ref{lm43}$ and
Corollaries~\ref{lm42}--\ref{lm429} reduce the Littlewood--Offord
problem to the study of the measure~$M^*$, uniquely corresponding
to the vector~$a$. In fact, almost all the results obtained when
solving this problem, are formulated in terms of the
coefficients~$a_j$ or, equivalently, in terms of the properties of
the measure~$M^*$. Sometimes this leads to a loss of information
on the distribution of the random variable~$X$, which can help in
obtaining more precise estimates. In particular, if $\mathcal
L(X)$ is the standard normal distribution, then $F_a$ is a
Gaussian distribution with zero mean and covariance operator which
can be easily calculated. Thus, there are situations in which it
is possible to obtain estimates for $ Q (F_a, \tau) $ that do not
follow from the results formulated in terms of the measure~$M^*$.

Note that using the results of Arak
 \cite{Arak:1980},~\cite{Arak:1981} (see also \cite{Arak:Zaitsev:1986})
one could derive from Theorem \ref{lm43} estimates similar to
estimates of concentration functions in the Littlewood--Offord
problem, which were obtained in a recent paper of Nguyen and
Vu~\cite{Nguyen:Vu:2011} (see also \cite{Nguyen and Vu13}). A
detailed discussion of this fact is presented in a joint paper of
the authors and Friedrich G\"otze which is preparing for the
publication. In the same paper there is a proof of
multidimensional analogs of some results of Arak \cite{Arak:1980}.
In Theorems~\ref{t4} and~\ref{thm1} below, we provide without
proof the formulations of these results which demonstrates a
relation between the order of smallness of the concentration
function of the sum and the arithmetic structure of the supports
of distributions of independent random vectors for {\it arbitrary}
distributions of summands, in contrast to the results of
\cite{Fried:Sod:2007}, \cite{Nguyen:Vu:2011},
\cite{Rud:Ver:2008}--\cite{Ver:2011}, in which a similar
relationship was found in a particular case of summands with the
distributions arising in the Littlewood--Offord problem.

We need some notation. Let ${\mathbf Z}_+$ be the set non-negative
integers. For any $r\in{\mathbf Z}_+$ and $u=(u_1,\ldots, u_r)\in
{({\mathbf R}^d)}^r$, $u_j\in {\mathbf R}^d$, $j=1,\ldots,r$,
introduce the sets
\begin{equation}
{K}_{1}(u)=\Big\{\sum\limits_{j=1}^r n_j u_j:n_j\in \{-1,0,1\}
\hbox{ for }j=1,\ldots,r\Big\}.\label{1s17}
\end{equation} We denote by $[B]_\tau$ the closed $\tau$-neighborhood of
a set $B$ in the sense of the norm $|\,\cdot\,|$.

\begin{theorem}\label{t4}
Let  $\tau\ge0$ and let $F_j$, $j=1,\ldots,n$, be  $d$-dimensional
probability distributions. Denote $\gamma=Q\Big(\prod_{j=1}^n
F_j,\tau\Big)$. Then there exist $r\in\mathbf Z_+$ and vectors
$u_1,\ldots, u_r;x_1,\ldots, x_r\in {{\mathbf R}^d}$ such that
\begin{equation}
r\ll_d \left|\log \rho\right|+1, \label{1s682}
\end{equation}
and
\begin{equation}
\sum\limits_{j=1}^n F_j\{{\mathbf R}^d\setminus[K_{1}(u)]_\tau+x_j\}\ll_d \bigl(\left|\log \rho\right|+1\bigr)^3,
\label{1s692}
\end{equation}
where $u=(u_1,\ldots, u_r)\in {({\mathbf R}^d)}^r$, and the set
$K_{1}(u)$ is defined in~\eqref{1s17}.
 \end{theorem}

\begin{theorem}\label{thm1}
Let $D$ be a $d$-dimensional infinitely divisible distribution
with characteristic function of the form
$\exp\big\{\alpha(\widehat M(t)-1)\big\}$, $t\in{\mathbf R}^d$,
where $\alpha>0$ and $M$ is a probability distribution. Let
$\tau\ge0$ and  $\gamma=Q(D, \tau)$. Then there exist $r\in\mathbf
Z_+$ and vectors $u_1,\ldots, u_r\in {{\mathbf R}^d}$ such that
\begin{equation}
r\ll_d \left|\log \gamma\right|+1, \label{1s65}
\end{equation}
and
\begin{equation}
\alpha \,M\{{\mathbf R}^d\setminus[K_{1}(u)]_\tau\}\ll_d
\bigl(\left|\log \gamma\right|+1\bigr)^3, \label{1s68d}
\end{equation}
where $u=(u_1,\ldots, u_r)\in {({\mathbf R}^d)}^r$.
\end{theorem}
\medskip

\noindent {\bf Proof of Theorem~\ref{lm43}.} Let us show that, for
arbitrary probability distribution~$F$ and $\lambda,T>0$,
\begin{multline}
\log\int_{|t|\le
T}\exp\Big(-\frac{\,1\,}2\;\sum_{k=1}^{n}\int_{z\in\mathbf{R}}\big(1-\cos(\left\langle
\, t,a_k\right\rangle z)\big)\,\lambda\,F\{dz\}\Big)\,dt\\ \le
\int_{z\in\mathbf{R}}\bigg(\log\int_{|t|\le
T}\exp\Big(-\frac{\,\lambda\,}2\;\sum_{k=1}^{n}\big(1-\cos(\left\langle
\, t,a_k\right\rangle z)\big)\Big)\,dt\bigg)\,F\{dz\}\\ =
\int_{z\in\mathbf{R}}\bigg(\log\int_{|t|\le
T}\widehat{H}_{z}^{\lambda}(t)\,dt\bigg)\,F\{dz\}.\label{dd11}
\end{multline}
It suffices to prove \eqref{dd11} for discrete  distributions~$F=
\sum_{j=1}^{\infty}p_j E_{z_j} $, where $0\le p_j\le1$,
$z_j\in\mathbf{R}$, $\sum_{j=1}^{\infty}p_j =1 $. Applying in this
case the H\"older inequality, we have
\begin{multline}
\int\limits_{|t|\le
T}\exp\Big(-\frac{\,1\,}2\;\sum\limits_{k=1}^{n}\int\limits_{z\in\mathbf{R}}\big(1-\cos(\left\langle
\, t,a_k\right\rangle z)\big)\,\lambda\,F\{dz\}\Big)\,dt\\ =
\int\limits_{|t|\le
T}\exp\Big(-\frac{\,\lambda\,}2\;\sum\limits_{j=1}^{\infty}p_j\sum\limits_{k=1}^{n}\big(1-\cos(\left\langle
\, t,a_k\right\rangle z_j)\big)\Big)\,dt\\ \le\prod\limits_{j=1}^{\infty}
\bigg(\int\limits_{|t|\le
T}\exp\Big(-\frac{\,\lambda\,}2\;\sum\limits_{k=1}^{n}\big(1-\cos(\left\langle
\, t,a_k\right\rangle z_j)\big)\Big)\,dt\bigg)^{p_j}.\label{d11}
\end{multline}
Taking the logarithms of the left and right-hand sides
of~\eqref{d11}, we get~\eqref{dd11}. In general case we can
approximate the distribution~$F$ by discrete distributions in the
sense of weak convergence and to pass to the limit. We use that
the weak convergence of probability distributions is equivalent to
the convergence of characteristic functions which is uniform on
bounded sets. Moreover, the weak convergence of symmetric
infinitely divisible distributions is equivalent to the weak
convergence of the corresponding spectral measure. Note also that
the integrals $\int_{|t|\le T}$ may be replaced in~\eqref{dd11} by
the integrals $\int_{t\in B}$ over an arbitrary Borel set~$B$.

 Since for characteristic function
$\widehat{W}(t)$ of a random vector $Y$, we have
$$|\widehat{W}(t)|^2 = \mathbf{E}\exp(i\langle
\,t,\widetilde{Y}\rangle ) = \mathbf{E}\cos(\langle
\,t,\widetilde{Y}\rangle ),$$ where $\widetilde{Y}$ is the
corresponding symmetrized random vector, then
\begin{equation}\label{6}|\widehat{W}(t)| \leq
\exp\Big(-\cfrac{\,1\,}{2}\,\big(1-|\widehat{W}(t)|^2\big)\Big)  =
\exp\Big(-\cfrac{\,1\,}{2}\,\mathbf{E}\,\big(1-\cos(\langle
\,t,\widetilde{Y}\rangle )\big)\Big).
\end{equation}

According to Theorem \ref{lm3} and relations $V=\lambda\,F\le G$,
\eqref{dd11} and \eqref{6}, we have
\begin{eqnarray}
Q(F_a,\tau)&\ll_{d}& \tau^d\int_{\tau|t|\le 1}|\widehat{F}_a(t)|\,dt \nonumber\\
&\ll_{d}& \tau^d\int_{\tau|t|\le
1}\exp\Big(-\frac{\,1\,}{2}\,\sum_{k=1}^{n}\mathbf{E}\,\big(1-\cos(\left\langle
\,t,a_k\right\rangle \widetilde{X})\big)\Big)\,dt\nonumber\\
&=&\tau^d\int_{\tau|t|\le
1}\exp\Big(-\frac{\,1\,}2\;\sum_{k=1}^{n}\int_{z\in\mathbf{R}}\big(1-\cos(\left\langle
\, t,a_k\right\rangle z)\big)\,G\{dz\}\Big)\,dt\nonumber\\
&\le&\tau^d\int_{\tau|t|\le
1}\exp\Big(-\frac{\,1\,}2\;\sum_{k=1}^{n}\int_{z\in\mathbf{R}}\big(1-\cos(\left\langle
\, t,a_k\right\rangle z)\big)\,\lambda\,F\{dz\}\Big)\,dt\nonumber\\
&\le&\exp\bigg(
\int_{z\in\mathbf{R}}\log\bigg(\tau^d\int_{\tau|t|\le
1}\widehat{H}_{z}^{\lambda}(t)\,dt\bigg)\,F\{dz\}\bigg).\label{cc11}
\end{eqnarray}

Using \eqref{8jrt} and  \eqref{1b}, we have
\begin{eqnarray}
 \tau^d\int_{\tau|t|\le 1}
\widehat{H}_{z}^{\lambda}(t) \,dt &\asymp_{d}&
Q(H_{z}^{\lambda},\tau) = Q\big(H_{1}^{\lambda},
\tau{|z|}^{-1}\big)\nonumber
\\&\le& \big(1+\lfloor\tau(\varepsilon |z|)^{-1}\rfloor \big)^d\, Q(H_{1}^{\lambda},\varepsilon).\label{cc22}
\end{eqnarray}
 Substituting this estimate into
\eqref{cc11}, we obtain~\eqref{cc23}. $\square$

\end{document}